\input amstex 
\documentstyle{amsppt} 
\loadbold
\magnification=1200

\pagewidth{6.43truein}
\pageheight{8.5truein}

\define\A{{\Cal A}}
\redefine\B{{\Cal B}}
\define\C{{\Bbb C}}
\redefine\D{{\Bbb D}}
\redefine\O{{\Cal O}}
\redefine\P{{\Bbb P}}
\define\R{{\Bbb R}}
\define\T{{\Bbb T}}
\redefine\phi{{\varphi}}
\redefine\epsilon{{\varepsilon}}
\define\ac{\acuteaccent}

\redefine\cdot{{\boldsymbol\cdot}}

\hyphenation{pluri-sub-har-mon-ic}

\refstyle{B}
\NoRunningHeads
\TagsOnRight

\topmatter 
\title The Siciak-Zahariuta Extremal Function \\ 
As the Envelope of Disc Functionals \endtitle
\author Finnur L\ac arusson and Ragnar Sigurdsson \endauthor
\address Department of Mathematics, University of Western Ontario,
London, Ontario N6A~5B7, Canada \endaddress
\email larusson\@uwo.ca \endemail
\address Science Institute, University of Iceland, Dunhaga 3,
IS-107 Reykjav\ac ik, Iceland \endaddress
\email ragnar\@hi.is \endemail 

\thanks The first-named author was supported in part by the Natural
Sciences and Engineering Research Council of Canada.  \endthanks

\thanks First version 22 April 2005.  Second, expanded version 6 July
2005.
\endthanks

\abstract We establish disc formulas for the Siciak-Zahariuta extremal
function of an arbitrary open subset of complex affine space.  This
function is also known as the pluricomplex Green function with
logarithmic growth or a logarithmic pole at infinity.  We extend
Lempert's formula for this function from the convex case to the
connected case.
\endabstract


\endtopmatter

\document

\specialhead Introduction \endspecialhead

\noindent The Siciak-Zahariuta extremal function $V_X$ of a subset $X$
of complex affine space $\C^n$ is defined as the supremum of all
entire plurisubharmonic functions $u$ of minimal growth with $u|X\leq
0$.  It is also called the pluricomplex Green function of $X$ with
logarithmic growth or a logarithmic pole at infinity (although this is
a bit of a misnomer if $X$ is not bounded).  A plurisubharmonic
function $u$ on $\C^n$ is said to have minimal growth (and belong to
the class $\Cal L$) if $u-\log^+\|\cdot\|$ is bounded above on $\C^n$.
If $X$ is open and nonempty, then $V_X\in \Cal L$.  More generally, if
$X$ is not pluripolar, then the upper semicontinuous regularization
$V_X^*$ of $V_X$ is in $\Cal L$, and if $X$ is pluripolar, then
$V_X^*=\infty$.  Siciak-Zahariuta extremal functions play a
fundamental role in pluripotential theory and have found important
applications in approximation theory, complex dynamics, and elsewhere.
For a detailed account of the basic theory, see \cite{K, Chapter 5}.
For an overview of some recent developments, see \cite{Pl}.

The extremal functions of pluripotential theory are usually defined as
suprema of classes of plurisubharmonic functions with appropriate
properties.  The theory of disc functionals, initiated by Poletsky in
the late 1980s \cite{P1, PS}, offers a different approach to extremal
functions, realizing them as envelopes of disc functionals.  A disc
functional on a complex manifold $Y$ is a map $H$ into
$\overline\R=[-\infty,\infty]$ from the set of analytic discs in $Y$,
that is, holomorphic maps from the open unit disc $\Bbb D$ into $Y$.
We usually restrict ourselves to analytic discs that extend
holomorphically to a neighbourhood of the closed unit disc.  The
envelope $EH$ of $H$ is the map $Y\to\overline\R$ that takes a point
$x\in Y$ to the infimum of the values $H(f)$ for all analytic discs
$f$ in $Y$ with $f(0)=x$.  Disc formulas have been proved for such
extremal functions as largest plurisubharmonic minorants, including
relative extremal functions, and pluricomplex Green functions of
various sorts, and used to establish properties of these functions
that had proved difficult to handle via the supremum definition.  Some
of this work has been devoted to extending to arbitrary complex
manifolds results that were first proved for domains in $\C^n$.  See
for instance \cite{BS, E, EP, LS1, LS2, LLS, P2, P3, R, RS}.

In the convex case, there is a disc formula for the Siciak-Zahariuta
extremal function due to Lempert \cite{M, Appendix}.  The main
motivation for the present work was to generalize Lempert's formula.
Because of the growth condition in the definition of the
Siciak-Zahariuta extremal function, we did not see how to fit it into
the theory of disc functionals until we realized, from a remark of
Guedj and Zeriahi \cite{GZ}, that minimal growth is nothing but
quasi-plurisubharmonicity with respect to the current of integration
along the hyperplane at infinity.  This observation is implicit in the
proof of our first main result, Theorem 1, which presents a family of
new disc formulas for the Siciak-Zahariuta extremal function of an
arbitrary open subset of affine space.  Theorem 2 contains more such
formulas.  Our second main result, Theorem 3, establishes Lempert's
formula, in slightly modified form, for every connected open subset of
affine space.  The formula is easily seen to fail for disconnected
sets in general.

Let us briefly summarize the contents of the paper.  We let $X$ be an
open subset of $\C^n$ and seek a disc formula for $V_X$.  If we have a
good upper semicontinuous majorant for $V_X$ on $\C^n$, so good that
$V_X$ is its largest plurisubharmonic minorant, then we have a disc
formula for $V_X$ as the so-called Poisson envelope of the majorant.
If $B$ is a ball in $X$, say the unit ball, then such a majorant is
easily seen to be given as zero on $X$ and $V_B=\log\|\cdot\|$ outside
$X$.  The first main idea is to introduce certain good sets of
analytic discs in complex projective space $\P^n$, adapted to $X$, and
get many more good majorants for $V_X$ as the envelopes of a new disc
functional (called $J$ below) over such sets.  The second main idea is
that the disc formulas for $V_X$ thus obtained are in fact closely
related to Lempert's formula in the convex case, even though they look
quite different at first sight.  The relationship appears when we
modify the Poisson functional by adding to it the non-negative
functional $J$, balancing this by taking the envelope over the larger
class of all analytic discs in $\P^n$.  We show that the envelope is
still $V_X$.  If we restrict to analytic discs in $\P^n$ with boundary
in $X$, the Poisson term disappears and $J$ alone remains.  This is
essentially Lempert's formula, so the envelope is still $V_X$ if $X$
is convex.  This is easily seen to fail in general if $X$ is
disconnected.  A proof of Lempert's formula in slightly modified form,
assuming only that $X$ is connected, concludes the paper.  The proof
relies on a judicious choice of a good set of analytic discs, as well
as a fundamental argument in the theory of disc functionals,
Poletsky's proof of the plurisubharmonicity of the Poisson envelope,
adapted here to a somewhat different purpose.

\specialhead Good sets of analytic discs and the first disc
formula \endspecialhead

\noindent
If $Y$ is a complex manifold, then we denote by $\A_Y$ the set of
analytic discs in $Y$, that is, the set of maps $\overline\D\to
Y$ that extend holomorphically to a neighbourhood of $\overline\D$.
If $H:\A_Y\to\overline\R$ is a disc functional on $Y$ and $\B\subset\A_Y$,
then the envelope of $H$ with respect to $\B$ is the function $E_\B H:
Y\to\overline\R$ with
$$E_\B H(y)=\inf\{H(f):f\in\B, f(0)=y\}, \qquad y\in Y.$$
We usually write $EH$ for $E_{\A_Y} H$ and simply call it the envelope
of $H$.

Perhaps the most important example of a disc functional is the Poisson
functional $f\mapsto \int_\T \phi\circ f\,d\sigma$ associated to an
upper semicontinuous function $\phi:Y\to[-\infty,\infty)$ (here,
$\sigma$ is the normalized arc length measure on the unit circle $\T$).
Its envelope is the largest plurisubharmonic minorant of $\phi$ on $Y$.
This was first proved for domains in affine space by Poletsky \cite{P1},
and later, with a different proof, by Bu and Schachermayer \cite{BS},
and finally generalized to all complex manifolds by Rosay \cite{R}.

We view $\C^n$ as the subset of $\P^n$ with projective
coordinates $[z_0:\dots:z_n]$ where $z_0\neq 0$ and write $H_\infty$
for the hyperplane at infinity where $z_0=0$.  We define a disc
functional $J$ on $\P^n$ by the formula
$$J(f)=-\sum_{\zeta\in f^{-1}(H_\infty)} m_{f_0}(\zeta)\log|\zeta|\ \geq
\ 0, \qquad f\in\A_{\P^n}.$$
Here, $m_{f_0}(\zeta)$ denotes the multiplicity of the intersection of
$f$ with $H_\infty$, that is, the order of the zero of the component
$f_0$ at $\zeta$ when $f$ is expressed as $[f_0:\dots:f_n]$ in
projective coordinates.  When the zeros of $f_0$ are not isolated, that
is, $f(\D)\subset H_\infty$, we take $J(f)=\infty$, and when
$f(\D)\cap H_\infty=\varnothing$, we take $J(f)=0$.

To indicate the relevance of $J$ to the Siciak-Zahariuta extremal
function, let $X\subset\C^n$ be open and $f\in\A_{\P^n}$ have
$f(\T)\subset X$.  For simplicity, we assume that $f$ sends only one
point $\zeta\in\D$ to $H_\infty$ and $m_{f_0}(\zeta)=1$.  Let $\rho$ be
the reciprocal and $\tau$ be an automorphism of $\D$ interchanging
$0$ and $\zeta$.  Then $g=f\circ\tau\circ\rho:\C\setminus\D\to\C^n$ is
holomorphic with a simple pole at infinity and $g(\T)\subset X$.  Hence,
$V_X\circ g$, extended as zero across $\D$, is a subharmonic function
on $\C$ of minimal growth, so $V_X\circ g\leq V_{\D}=\log|\cdot|$ and
$$V_X(f(0)) = V_X(g(1/\zeta)) \leq -\log|\zeta| = J(f).$$

A subset $\B$ of $\A_{\P^n}$ is called {\it good} with respect to
an open subset $X$ of $\C^n$ if:
\roster
\item $f(\T)\subset X$ for every $f\in\B$,
\item for every $z\in\C^n$, there is a disc in $\B$ with centre $z$,
\item for every $x\in X$, the constant disc at $x$ is in $\B$, and
\item the envelope $E_\B J$ is upper semicontinuous on $\C^n$ and has
minimal growth, that is, $E_\B J-\log^+\|\cdot\|$ is bounded above on
$\C^n$.
\endroster
Note that by (2), $0\leq E_\B J<\infty$ on $\C^n$, by (3), $E_\B J=0$
on $X$, and clearly $E_\B J=\infty$ on $H_\infty$.  Property (4) may be
hard to verify directly, but Proposition 2 gives a useful sufficient
condition for it to hold.  Roughly speaking, if $\B$ contains a disc
centred at each point of $\P^n$, then (4) holds if (but not only if)
discs in $\B$ can be varied continuously.

\proclaim{Theorem 1}  Let $X$ be an open subset of $\C^n$ and $\B$ be a
good set of analytic discs in $\P^n$ with respect to $X$.  Then
the Siciak-Zahariuta extremal function $V_X$ of $X$ is the envelope of
the disc functional $H_\B$ on $\P^n$ defined by the formula
$$H_\B(f)=J(f)+\int_{\T\setminus f^{-1}(X)}E_\B J\circ f \,d\sigma,
\qquad f\in\A_{\P^n}.$$
\endproclaim

\noindent
{\bf Remarks. 1.}  We define $V_X=\infty$ on $H_\infty$
and it is clear that $EH_\B=\infty$ on $H_\infty$.  Since $E_\B J=0$
on $X$, the integral above might as well be taken over all of $\T$.
The disc functional $H_\B$ is thus given as the Poisson functional of
$E_\B J$ minus the Lelong-like functional $-J$ (see \cite{LS2}
for the definition of the Lelong functional).  Envelopes of disc
functionals associated to complex subspaces in the way that $-J$ is
associated to $H_\infty$ are Green functions of a type studied in
\cite{RS}.

{\bf 2.}  Using Proposition 2, it is easy to see that
the largest subset $\B\subset\A_{\P^n}$ that is good with respect
to $X$ is the set $\A_{\P^n}^X$ of all $f\in\A_{\P^n}$ with
$f(\T)\subset X$.  This yields the smallest possible $E_\B J$ in the
second term of the formula for $H_\B$.  We write $H_X$ for
$H_{\A_{\P^n}^X}$, so $H_X$ is the smallest disc functional $H_\B$ where
$\B\subset\A_{\P^n}$ is good with respect to $X$.  Other choices of $\B$
make the second term explicitly computable and yield information about
almost extremal discs (see Propositions 4 and 5).  Theorem 2 shows that
$V_X$ is in fact the envelope of $H_\B$ over analytic discs in
$\C^n$ only; for such discs, $J$ vanishes.

{\bf 3.} By a result of Lempert \cite{M, Appendix}, if $X$ is convex,
then $V_X$ is the envelope of $H_\B$ over analytic discs in $\P^n$
with boundary in $X$ and at most one simple pole; for such discs, the
second term vanishes, leaving only $J$.  We discuss this in detail
later in the paper.  For disconnected $X$, it is generally not true
that $V_X=E_{\A_{\P^n}^X}J$.  For example, suppose $X$ is the disjoint
union of two nonempty convex open sets $Y$ and $Z$.  Then
$$E_{\A_{\P^n}^X}J = \min\{E_{\A_{\P^n}^Y}J, E_{\A_{\P^n}^Z}J\}
=\min\{V_Y, V_Z\}$$
is not even plurisubharmonic in general (but it does provide an
upper bound for $V_X$).

\demo{Proof of Theorem 1}  Let $\pi:Z=\C^{n+1}\setminus\{0\}\to\P^n$ be
the projection.  Write $Z_0=\pi^{-1}(H_\infty)=\{z\in Z:z_0=0\}$.
The advantage of working on $Z$ rather than on $\P^n$
is that the pullback of the current of integration along $H_\infty$ has
a global plurisubharmonic potential $\phi(z)=\log|z_0|$ on $Z$.  Note that
if $x\in\P^n$ and $z\in Z$ with $\pi(z)=x$, then every analytic disc in
$\P^n$ centred at $x$ lifts to an analytic disc in $Z$ centred at $z$.
Hence, as $f$ runs through all analytic discs in $Z$ with $f(0)=z$,
$\pi\circ f$ runs through all analytic discs $g$ in $\P^n$ with $g(0)=x$.  

Let $\pi^*\B$ be the set of analytic discs $f$ in $Z$ with $\pi\circ
f\in\B$.  Define a function $\psi$ on $Z\setminus Z_0$ by the formula
$$\psi(z)=\inf\{\int_\T \phi\circ f\,d\sigma: f\in\pi^*\B,
f(0)=z\}, \qquad z\in Z\setminus Z_0.$$
By the defining property (2) of a good set of analytic discs and
plurisubharmonicity of $\phi$, we have $\phi\leq\psi<\infty$ on
$Z\setminus Z_0$, and by property (3), $\psi=\phi$ on $\pi^{-1}(X)$.

If $f\in\pi^*\B$, $f(0)\not\in Z_0$, then the Riesz Representation
Theorem applied to the subharmonic function $\phi\circ f=\log|f_0|$
gives
$$\phi(f(0)) = \int_\T\phi\circ f\,d\sigma + \tfrac
1{2\pi}\int_\D \log|\,\cdot\,|\Delta(\phi\circ f).$$
Also,
$$\tfrac 1{2\pi}\int_\D \log|\,\cdot\,|\Delta(\phi\circ f)
= \sum_{\zeta\in f_0^{-1}(0)} m_{f_0}(\zeta)\log|\zeta|
=-J(\pi\circ f),$$
so
$$\int_\T\phi\circ f\,d\sigma=\phi(f(0))+J(\pi\circ f),$$
and
$$\psi=\phi + E_\B J \circ\pi \qquad\text{on }Z\setminus Z_0.$$
By property (4), $E_\B J$ is upper semicontinuous on
$\C^n$, so $\psi:Z\setminus Z_0 \to\R$ is upper semicontinuous.  The
minimal growth condition on $E_\B J$ means that $E_\B J\circ\pi+\phi
= \psi$ is locally bounded above at $Z_0$.  Hence, the upper
semicontinuous extension $\psi^*:Z\to[-\infty,\infty)$, which we shall
simply call $\psi$, is well defined, and we have $\phi\leq\psi$ on $Z$
and $\psi=\phi$ on $\pi^{-1}(X)$.  The key property of $\psi$ is that
if $u$ is a plurisubharmonic function on $Z$ and $u\leq\phi$ on
$\pi^{-1}(X)$, then $u\leq\psi$ on $Z\setminus Z_0$ by property
(1), so $u=(u|Z\setminus Z_0)^*\leq\psi$ on all of $Z$.
The converse is clear since $\psi=\phi$ on $\pi^{-1}(X)$.

Now $u\in\Cal L$ if and only if $u\circ\pi+\phi$ is plurisubharmonic on
$Z$.  Namely, the minimal growth condition that defines $\Cal L$ means
that $u\circ\pi+\phi$, which is plurisubharmonic on $Z\setminus Z_0$, is
locally bounded above at $Z_0$, which in turn means that $u\circ\pi+\phi$
extends uniquely to a plurisubharmonic function on all of $Z$.

Hence, $u\in\Cal L$ and $u\leq 0$ on $X$ if and only if $u\circ\pi
+\phi$ is plurisubharmonic on $Z$ and $u\circ\pi+\phi\leq\phi$ on
$\pi^{-1}(X)$, that is, $u\circ\pi+\phi\leq\psi$ on $Z$.  Thus, clearly,
$V_X\circ\pi+\phi\leq\psi$.  Also, since $\psi-\phi=E_\B J\circ\pi$ on
$Z\setminus Z_0$ is invariant under homotheties, so is its largest
plurisubharmonic minorant $P_{Z\setminus Z_0}(\psi-\phi)$ on $Z\setminus
Z_0$.  Since $\phi$ is pluriharmonic on $Z\setminus Z_0$, $P_{Z\setminus
Z_0}(\psi-\phi) = P_{Z\setminus Z_0}(\psi)-\phi$, so $P_{Z\setminus Z_0}
(\psi)=u\circ\pi+\phi$, where $u\in\Cal L$ and $u\leq 0$ on $X$.
Therefore, $u\leq V_X$ and $P_Z(\psi)|Z\setminus Z_0 \leq P_{Z\setminus
Z_0}(\psi)\leq V_X\circ\pi+\phi$ on $Z\setminus Z_0$, so $P_Z(\psi)\leq
V_X\circ\pi+\phi$ on $Z$.

This shows that $V_X\circ\pi+\phi$ is the largest plurisubharmonic
minorant of $\psi$ on $Z$, so Poletsky's theorem yields a disc
formula for $V_X\circ\pi+\phi$ as the Poisson envelope of $\psi$ on $Z$.
For $z\in Z\setminus Z_0$, it follows that $V_X(\pi(z))$ is the infimum
over all analytic discs $f$ in $Z$ with $f(0)=z$ of the numbers
$$\int_\T\psi\circ f\,d\sigma - \phi(z).$$
By the Riesz Representation Theorem,
$$\align \int_\T\psi\circ f\,d\sigma - \phi(z) &= \int_\T(\psi-\phi)
\circ f\,d\sigma + J(\pi\circ f) \\ &= J(\pi\circ f)+\int\limits_{\T
\setminus f^{-1}(\pi^{-1}(X))} E_\B J\circ\pi\circ f\,d\sigma.
\endalign$$
Note that $f^{-1}(Z_0)\cap\T$ is finite, so the second and third
integrals are equal.  This shows that $V_X(x)=EH_\B(x)$ for all
$x\in\C^n$.  For $x\in H_\infty$, this is obvious, as mentioned in
Remark 1 above.
\qed\enddemo

The multiplicity factor in the definition of $J$ may be omitted without
affecting Theorem 1 with $\B=\A_{\P^n}^X$, that is, without changing
$E_{\A_{\P^n}^X} J$ or $EH_X$.

\proclaim{Proposition 1}  Let $X\subset\C^n$ be open and $f\in \A_{\P^n}$
have $f(0)\not\in H_\infty$.  Then there is $g\in \A_{\P^n}$ with
$g(0)=f(0)$, $m_{g_0}=1$ on $g^{-1}(H_\infty)$, and $J(f)=J(g)$, such
that $g$ is uniformly as close to $f$ on $\overline\D$ as we wish, so
in particular, if $f(\T)\subset X$, then $g(\T)\subset X$.
\endproclaim

\demo{Proof}  Now $f$ intersects $H_\infty$ in finitely many
points $a_1,\dots,a_k\in\D\setminus\{0\}$ with multiplicities
$m_j=m_{f_0}(a_j)$.  Let $\tilde f\in \A_Z$ be a lifting
of $f$.  By exactly the same argument as in the proof of Lemma 3.1 in
\cite{LS2}, taking the function $\alpha$ there to be the characteristic
function of $Z_0$ in $Z$, we obtain $\tilde g\in \A_Z$ arbitrarily
uniformly close to $\tilde f$ on $\overline\D$ such that $\tilde g(0)
=\tilde f(0)$, the zeros $c_1,\dots,c_m$ of $\tilde g_0$ in $\D$ all
have multiplicity $1$, their number $m$ equals $m_1+\cdots+m_k$, and
$$\sum_{j=1}^m\log|c_j|=\sum_{j=1}^k m_j\log|a_j|.$$
Finally, take $g=\pi\circ \tilde g$.
\qed\enddemo

\specialhead  Further results on good sets of analytic discs
\endspecialhead

\noindent
Using the proof of Theorem 1, we now present a sufficient condition for
a set of analytic discs to be good.

\proclaim{Proposition 2}  Let $X$ be an open subset of $\C^n$ and
$\B$ be a subset of $\A_{\P^n}$ satisfying the following two properties:
\roster
\item"(2')" For every $z\in\P^n$, there is a disc in $\B$ with centre
$z$.
\item"(4')" Discs in $\B$ can be varied continuously, that is, for every
$f\in\B$ there is a map from a neighbourhood $U$ of $f(0)$ into
$\B$, continuous as a map $\overline\D\times U \to\P^n$, taking each
$x\in U$ to a disc centred at $x$ and taking $f(0)$ to $f$.
\endroster
Then $\B$ satisfies property {\rm (4)} in the definition of a good set
of analytic discs.

Hence, if $\B$ satisfies {\rm (1)}, {\rm (2')}, {\rm (3)}, and
{\rm (4')}, then $\B$ is good with respect to $X$.
\endproclaim

\demo{Proof}  Define
$$\psi(z)=\inf\{\int_\T \phi\circ f\,d\sigma: f\in\pi^*\B, f(0)=z\},
\qquad z\in Z\setminus Z_0,$$
as in the proof of Theorem 1.  Properties (2) and (4') imply that
$\psi:Z\setminus Z_0\to\R$ is upper semicontinuous.  Now
$E_\B J \circ\pi=\psi-\phi$ on $Z\setminus Z_0$, so $E_\B J$ is upper
semicontinuous on $\C^n$.  Moreover, $E_\B J$ has minimal growth since
$\psi$ is locally bounded above at $Z_0$ by (2') and (4').
\qed\enddemo

The next result provides an interesting class of examples of good sets
of analytic discs.

\proclaim{Proposition 3}  Let $X$ be a connected open subset of $\C^n$
and let $\beta$ be a free homotopy class of loops in $X$, that is, of
continuous maps $\T\to X$.  Let $\B$ be the set of analytic discs $f$
in $\P^n$ such that $f(\T)\subset X$ and $f|\T\in\beta$.  Then $\B$
satisfies properties {\rm (1)}, {\rm (2')}, and {\rm (4')}.  If $\beta$
is the trivial class, then $\B$ also satisfies property {\rm (3)}, so
$\B$ is good with respect to $X$.
\endproclaim

\demo{Proof}  Only (2') is not obvious.  Let $z\in\C^n$ and a continuous
map $\alpha:\T\to X$ be a representative for $\beta$.  Rational functions
on $\C$ whose poles lie outside $\T\cup\{0\}$ are uniformly dense among
continuous functions $\T\cup\{0\}\to\C$ (see e.g\. \cite{AW, Theorem 2.8}).
Therefore, for each $\epsilon>0$, we obtain rational functions
$f_1,\dots,f_n$ without poles on $\T\cup\{0\}$, defining an analytic
disc $f=(f_1,\dots,f_n)$ in $\P^n$, such that $f(0)=z$ and $f|\T$ is
within $\epsilon$ of $\alpha$, so $f|\T$ is freely homotopic to $\alpha$
in $X$ if $\epsilon$ is small enough.  If $z\in H_\infty$, we reduce to
the previous case by moving $z$ into $\C^n$ by an automorphism of $\P^n$
close to the identity.
\qed\enddemo

\specialhead  Majorants for the Siciak-Zahariuta function and the second
disc formula \endspecialhead

\noindent
Let $X$ be an open subset of $\C^n$ and $\B$ be a good set of analytic
discs in $\P^n$ with respect to $X$.  By Theorem 1, $V_X=EH_\B$.
Clearly, $EH_\B\leq E_\B H_\B=E_\B J$, so $V_X\leq E_\B J$.  Moreover,
if $u$ is a plurisubharmonic function on $\C^n$ with $u\leq E_\B J$,
then $u\leq 0$ on $X$ by property (3) in the definition of a good set of
analytic discs, and $u$ has minimal growth by property (4), so $u\leq
V_X$.  This shows that $V_X$ is the largest plurisubharmonic minorant,
and hence the Poisson envelope, of $E_\B J$ on $\C^n$.  It follows that
$E_\B J$ is plurisubharmonic if and only if $E_\B J = V_X$.  We have
proved the following result.

\proclaim{Theorem 2}  Let $X$ be an open subset of $\C^n$ and $\B$ be a
good set of analytic discs in $\P^n$ with respect to $X$.  Then $V_X$ is
the largest plurisubharmonic minorant of $E_\B J$ on $\C^n$.
Consequently, for every $z\in \C^n$,
$$V_X(z)=\inf \int_{\T\setminus f^{-1}(X)}E_\B J\circ f \,d\sigma,$$
where the infimum is taken over all analytic discs $f$ in $\C^n$ with
$f(0)=z$.
\endproclaim

\specialhead The third disc formula and almost extremal discs
\endspecialhead

\noindent
Let $X$ be an open subset of $\C^n$.  The simple disc formula for $V_X$
mentioned in the Introduction is in fact a special case of Theorem 2.
Namely, suppose $B$ is a closed ball with centre $a$ and radius $r>0$
contained in $X$.  As is well known, $V_B=\log\|\cdot-a\|-\log r$
outside $B$.  Setting $w=V_B$ outside $X$ and $w=0$ on $X$, we
obtain an upper semicontinuous majorant $w:\C^n\to[0,\infty)$ for
$V_X$.  It is easily seen that $V_X$ is the largest plurisubharmonic
minorant and hence the Poisson envelope of $w$.  Let us record this
fact.

\proclaim{Proposition 4}  Let $X$ be an open subset of $\C^n$ containing
the closed ball with centre $a$ and radius $r>0$.  Then $V_X$ is the
envelope of the disc functional $H_r$ on $\C^n$ defined by the formula
$$H_r(f)= \int_{\T\setminus f^{-1}(X)} \log\|f-a\| \,d\sigma -
\sigma(\T\setminus f^{-1}(X))\log r, \qquad f\in\A_{\C^n}.$$
\endproclaim

For simplicity, let us assume that $a$ is the origin.  Let $\B$ contain
the constant analytic discs in $X$ as well as the analytic discs $g_z$ in
$\P^n$ with
$$g_z(\zeta)=\frac{\|z\|+r\zeta}{r+\|z\|\zeta}\frac{r}{\|z\|}\, z$$
for each $z\in\C^n\setminus X$.  Note that
$g_z$ is centred at $z$, lies in the projective line through $z$
and the origin, and has its boundary on the sphere of radius $r$ centred
at the origin.  Also, $g_z$ sends one point in $\D$ to $H_\infty$, namely
$-r/\|z\|$.  Hence, $E_\B J(z)=J(g_z)=\log\|z\|-\log r$ if $z\in\C^n
\setminus X$, and $E_\B J=0$ on $X$, so $E_\B J=w$.  The defining
conditions for $\B$ to be a good set of analytic discs are easily
verified.  This shows that Proposition 4 is a special case of Theorem
2.  Note that the good set $\B$ satisfies neither property (2') nor (4')
in Proposition 2.

By the disconnected example in Remark 3 above, the following description
of almost extremal discs may be said to be optimal.  Namely, we
cannot always obtain $V_X$ as the envelope of $H_X$ or, equivalently,
of $J$ over analytic discs in $\P^n$, let alone $\C^n$, with boundaries
in $X$.  (Recall that $H_X$ was introduced as shorthand for
$H_{\A_{\P^n}^X}$ in Remark 2.)

\proclaim{Proposition 5}  Let $X$ be a nonempty open subset of $\C^n$.
Let $K$ be a compact subset of $X$ and $z\in\C^n$.  For each $\epsilon>0$,
there is an analytic disc $f$ in $\C^n$ centred at $z$ such that
$$V_X(z) \leq H_X(f) = \int_{\T\setminus f^{-1}(X)}E_{\A_{\P^n}^X}
J\circ f \,d\sigma < V_X(z)+\epsilon$$
and
$$\sigma(\T\setminus f^{-1}(X\setminus K))<\epsilon.$$
\endproclaim

\demo{Proof}  Say $X\setminus K$ contains a closed ball of radius $R>0$.
By Proposition 4 applied to $X\setminus K$, for each $0<r\leq R$, there
is $f_r\in\A_{\C^n}$ with $f_r(0)=z$ and
$$V_X(z) \leq H_X(f_r) \leq H_{X\setminus K}(f_r)
\leq H_r(f_r)<V_{X\setminus K}(z)+\epsilon=V_X(z)+\epsilon.$$  
Thus, as $r\to 0$, we must have $\sigma(\T\setminus f_r^{-1}
(X\setminus K))\to 0$, so we take $f=f_r$ with $r$ small enough.
\qed\enddemo

Since $V_X \leq E_{\A_{\P^n}^X} J$, Proposition 5 has the curious
consequence, for every open subset $X$ of $\C^n$, that $V_X$ is its own
Poisson envelope with respect to analytic discs in $\C^n$ that take all
but an arbitrarily small piece of the circle $\T$ into $X$.

\specialhead Relationship to the work of Lempert in the convex
case \endspecialhead

\noindent
We will now describe the relationship between Lempert's disc formula
for the Siciak-Zahariuta extremal function in the convex case, an
account of which was provided by Momm in \cite{M, Appendix}, and our
first disc formula.

Let $K$ be a strictly convex compact subset of $\C^n$ with real
analytic boundary and let $z\in\C^n\setminus K$.  Lempert's formula
states that $V_K(z)$ is the infimum of the numbers $\log r$ over all
holomorphic maps $f:\C\setminus\overline\D \to \C^n$ with a continuous
extension to $\T$ such that $f(\T)\subset K$, $f(r)=z$ with $r>1$, and
$\|f\|/|\cdot|$ is bounded, meaning that $f$ has at most a simple pole
at $\infty$.  (Furthermore, extremal maps exist and can be described
explicitly.)  Precomposing
$f$ with the reciprocal, we see that $V_K(z)$ is the infimum of the
numbers $-\log|\zeta|$ over all $f\in\A_{\P^n}$ with $f(\T)\subset K$
and $f(\zeta)=z$ such that $f$ maps into $\C^n$ except for at most a
simple pole at $0$.  Precomposing $f$ by an automorphism of $\D$ that
interchanges $0$ and $\zeta$, we see that $V_K(z)$ is the infimum of
the numbers $-\log|\zeta|$ over all $f\in\A_{\P^n}$ with $f(\T)\subset
K$ and $f(0)=z$ such that $f$ maps into $\C^n$ except for at most a
simple pole at $\zeta$.  For such a map $f$, we have
$-\log|\zeta|=J(f)$.

Let $X$ be a convex open subset of $\C^n$.  Then $X$ can be written as
the increasing union of relatively compact open subsets $X_n$, $n\geq
1$, such that the closure $\overline X_n$ is strictly convex with real
analytic boundary.  Namely, take a strictly convex exhaustion function
of $X$, such as the sum of $\|\cdot\|^2$ and the reciprocal of the
Euclidean distance to the boundary, and Weierstrass-approximate it by
a polynomial; the generic sublevel sets of the polynomial will be
smooth.

Now $V_X$ is the decreasing limit of $V_{X_n}$ and hence also the
decreasing limit of $V_{\overline X_n}$ as $n\to\infty$.  Therefore,
by Lempert's formula, $V_X(z)$ for $z\in\C^n\setminus X$, and thus
obviously for all $z\in\P^n$, is the infimum of the numbers $J(f)$
over all $f\in\A_{\P^n}$ with $f(\T)\subset X$ and $f(0)=z$ such that
$f$ maps into $\C^n$ except for at most one simple pole.  By Theorem
1, in between this infimum and $V_X(z)$ is the infimum of $J(f)$ over
the larger class of $f\in\A_{\P^n}$ with $f(\T)\subset X$ and
$f(0)=z$.

Thus, Lempert's formula can be stated as the following strengthening
of Theorem 1 for the convex case.

\proclaim{Lempert's formula} Let $X$ be a convex open subset of
$\C^n$.  Then the Siciak-Zahariuta extremal function of $X$ is the
envelope of $J$ with respect to the set of analytic discs in $\P^n$
with boundary in $X$ and at most one simple pole.  It follows that
$$V_X=E_{\A_{\P^n}^X}J.$$
\endproclaim

We conclude the paper by proving Lempert's formula for any connected
open subset $X$ of $\C^n$ in the slightly weakened form $V_X=
E_{\A_{\P^n}^X}J$, which we will henceforth refer to as Lempert's
formula.  As remarked earlier, this formula is easily seen to fail for
disconnected sets in general.  We do not know whether the stronger,
original form of Lempert's formula, using only analytic discs with at
most one simple pole, extends to all connected sets.

\specialhead Lempert's formula for arbitrary domains \endspecialhead

\noindent Let $X$ be a connected open subset of $\C^n$.  We may assume
that $X$ is neither empty nor all of $\C^n$ (otherwise, Lempert's
formula is obvious).  From now on, $\B$ will denote the set of
analytic discs in $\P^n$ containing all the constant discs in $X$ and
every disc
$$f_{z,w,r}:\zeta\mapsto w + \frac{\|z-w\|+r\zeta}{r+\|z-w\|\zeta}\,\frac
r{\|z-w\|} (z-w)$$ 
in $\P^n$, where $z\in\C^n\setminus X$, $w\in X$, and $r$ is less than
the Euclidean distance $d(w,\partial X)$ from $w$ to the boundary
$\partial X$ of $X$.  Observe that $f_{z,w,r}$ is injective, centred at
$z$, takes one point to $H_\infty$, namely $-r/\|z-w\|$, lies in the
projective line through $z$ and $w$, and maps $\T$ onto the circle
with centre $w$ and radius $r$ in this line.  It is easy to verify
that
$$E_\B J = \inf_{w\in X} \log^+\frac{\|\cdot-w\|}{d(w,\partial X)}
= \inf \{V_B : B \text{ is a ball in }X\}. $$
It follows that $\B$ is a good set of analytic discs in $\P^n$ with
respect to $X$.  Note that $E_\B J$ is not plurisubharmonic in
general: just consider an annulus.

If $X$ is smoothly bounded, then, using balls touching the boundary
from the inside, we see that $E_\B J=0$ on $\partial X$.  Now $V_X\leq
E_{\A_{\P^n}^X}J \leq E_\B J$ on $\C^n$, so $V_X = E_{\A_{\P^n}^X}J$
on $\overline X$.  Since every domain can be exhausted by smoothly
bounded domains and Lempert's formula is preserved by increasing
unions, it suffices to prove the formula on $\C^n\setminus\overline X$
assuming $E_\B J=0$ on $\partial X$.  The argument is based on the
following result.

\proclaim{Lemma} Let $X$ be a connected open subset of $\C^n$ and $\B$
be as above.  For every analytic disc $h$ in $\C^n\setminus\overline
X$, continuous function $v\geq E_\B J$ on $\C^n\setminus\overline X$, 
and $\epsilon>0$, there is $g\in\A_{\P^n}^X$ with $g(0)=h(0)$ and
$$J(g)\leq\int_\T v\circ h\,d\sigma+\epsilon.$$
\endproclaim

Fixing $z\in\C^n\setminus\overline X$ and taking the infimum over all
$v$, $\epsilon$, and $h$ with $h(0)=z$ as in the Lemma, we see that
$E_{\A_{\P^n}^X}J$ is no larger than the Poisson envelope, that is,
the largest plurisubharmonic minorant $P_{\C^n\setminus\overline
X}E_\B J$, of $E_\B J$ on $\C^n\setminus\overline X$.  Now
$$P_{\C^n\setminus\overline X}E_\B J = P_{\C^n}E_\B J
|\C^n\setminus\overline X.$$
Namely, if $u$ is plurisubharmonic on $\C^n\setminus\overline X$ and
$u\leq E_\B J$, then, after replacing $u$ by $\max\{u,0\}$ and using
the assumption that $E_\B J=0$ on $\partial X$, we can extend $u$ to a
plurisubharmonic function on all of $\C^n$ by setting $u=0$ on
$\overline X$.  Then $u\leq E_\B J$ on $\C^n$, so $u\leq P_{\C^n}E_\B
J$.  This proves one inequality; the other is obvious.  Finally, by
the remarks preceding Theorem 2, $P_{\C^n}E_\B J=V_X$ since $\B$ is good.
Thus, given the Lemma, we have established Lempert's formula:

\proclaim{Theorem 3} The Siciak-Zahariuta extremal function of a
connected open subset $X$ of $\C^n$ is the envelope of $J$ with
respect to the set of analytic discs in $\P^n$ with boundary in $X$,
that is,
$$V_X=E_{\A_{\P^n}^X}J.$$
\endproclaim

It remains to prove the Lemma.  Our argument is an adaptation of
Poletsky's orginal proof of the plurisubharmonicity of the Poisson
envelope.  See \cite{P1} or \cite{LS1, Section 2}.  We proceed as if
we were trying to show that $E_\B J$ was plurisubharmonic.

\demo{Proof of the Lemma} 
Take $\zeta_0\in \T$ and set $z_0=h(\zeta_0)$.  By the
definition of $\B$, there exist $w_0\in X$ and $r_0<d(w_0,{\partial}X)$ with
$$
J(f_{z_0,w_0,r_0})=\log(\|z_0-w_0\|/r_0)< E_\B J(z_0)+\epsilon.
$$ 
By continuity, there exists an open arc $I_0$ containing $\zeta_0$ such that
$$
J(f_{h(\zeta),w_0,r_0})=\log(\|h(\zeta)-w_0\|/r_0)<
v(h(\zeta))+\epsilon/2,
\qquad \zeta\in I_0.
$$ 
By compactness, there exist a cover of $\T$ by open arcs
$I_1,\dots,I_m$, points $w_1,\dots,w_m$ in $X$, and
$r_1,\dots,r_m>0$ such that $r_j<d(w_j,{\partial}X)$ and
$$
J(f_{h(\zeta),w_j,r_j})=\log(\|h(\zeta)-w_j\|/r_j)<
v(h(\zeta))+\epsilon/2, \qquad \zeta\in I_j, \ j=1,\dots,m.
$$ 
There exist $A\subset \{1,\dots,m\}$ and closed arcs $J_j\subset
I_j$, $j\in A$, which cover $\T$ with disjoint interiors.  By
possibly renumbering the arcs and splitting the interval $I_j$
containing $1$, we may assume that $A=\{1,\dots,m\}$ and
$$
J_j=\{e^{i\theta} : \theta\in [a_j,a_{j+1}]\}, \quad \text{ where }
\quad 0=a_1<a_2<\cdots<a_{m+1}=2\pi.
$$
Then
$$
\sum_{j=1}^m \int_{J_j} J(f_{h(\zeta),w_j,r_j})\, d\sigma(\zeta) <
\int_\T v\circ h \, d\sigma+\epsilon/2.
\tag{1}
$$ 
Since $X$ is connected, we can join $w_j$ and $w_{j+1}$ by a
$C^\infty$ path $\alpha_j:[0,1]\to X$ with $\alpha_j(0)=w_j$,
$\alpha_j(1)=w_{j+1}$, and choose a $C^\infty$ function
$\beta_j:[0,1]\to (0,\infty)$ with $\beta_j(0)=r_j$,
$\beta_j(1)=r_{j+1}$, and $\beta_j<d(\alpha_j,\partial X)$.
Here we take $w_{m+1}=w_1$ and $r_{m+1}=r_1$.  We may assume that the
derivatives of all orders of $\alpha_j$ and $\beta_j$ vanish at $0$ 
and $1$.  We choose
$$
C>\sum_{j=1}^m \sup_{\zeta\in J_j, t\in [0,1]} 
\big|J(f_{h(\zeta),w_j,r_j})-J(f_{h(\zeta),\alpha_j(t),\beta_j(t)}) 
\big| \tag{2}
$$ 
and $\delta>0$ such that $C\delta<\epsilon/2$ and
$\delta<\min_j(a_{j+1}-a_j)$.  We split each arc $J_j$ into the subarcs
$K_j=\{e^{i\theta} : \theta\in [a_j,a_{j+1}-\delta]\}$ and
$L_j=\{e^{i\theta} : \theta\in [a_{j+1}-\delta,a_{j+1}]\}$, and define
the $C^\infty$ loop $\gamma:\T\to X$ by
$$
\gamma(\zeta)=\cases w_j, &\zeta\in K_j,\  j=1,\dots,m,\\
\alpha_j((\theta-a_{j+1}+\delta)/\delta), &\zeta=e^{i\theta}\in L_j,
\ j=1,\dots,m,
\endcases
$$ 
the $C^\infty$ function $\varrho:\T\to (0,\infty)$ by
$$
\varrho(\zeta)=\cases r_j, &\zeta\in K_j, \  j=1,\dots,m,\\
\beta_j((\theta-a_{j+1}+\delta)/\delta), &\zeta=e^{i\theta}\in L_j,
\ j=1,\dots,m,
\endcases
$$ 
and, finally, the $C^\infty$ family 
$$
F(\cdot,\zeta)=f_{h(\zeta),\gamma(\zeta),\varrho(\zeta)}, \qquad
\zeta\in \T,
$$
of analytic discs in $\A_{\P^n}^X$.  By (1) and (2), 
$$
\int_\T J(F(\cdot,\zeta))\, d\sigma(\zeta) 
<\sum_{j=1}^m \int_{J_j} J(f_{h(\zeta),w_j,r_j})\,
d\sigma(\zeta)+C\delta
<\int_\T v\circ h \, d\sigma+\epsilon.
\tag{3}
$$ 
We take the lifting $\tilde h=(1,h)\in \A_Z$ of $h$ to
$Z=\C^{n+1}\setminus\{0\}$ by the projection $\pi:Z\to\P^n$, and the
lifting $\tilde f_{z,w,r}$ of $f_{z,w,r}$ given by
$$
\tilde f_{z,w,r}(\xi)=(\|z-w\|\xi/r+1,
(\|z-w\|\xi/r+1)w+(r\xi/\|z-w\|+1)(z-w)).
$$
Then the lifting $\tilde F(\cdot,\zeta)=\tilde
f_{h(\zeta),\gamma(\zeta),\varrho(\zeta)}$ of $F$ satisfies
$\tilde F(0,\cdot)=\tilde h$ on $\T$.

Take $r>1$ such that $h\in \O(D_r,\C^n)$ and $F(\cdot,\zeta)\in
\O(D_r,\P^n)$ for all $\zeta\in\T$, where $D_r=\{z\in\C:|z|<r\}$, and
define $\tilde F_j\in \O(D_r\times (D_r\setminus\{0\}), \C^{n+1})$,
$j\geq 1$, by
$$
\tilde F_j(\xi,\zeta)=\tilde h(\zeta)+\sum_{k=-j}^j
\bigg(\dfrac 1{2\pi}\int_0^{2{\pi}}
\big(\tilde F(\xi,e^{i{\theta}})-\tilde h(e^{i{\theta}})\big)
e^{-ik{\theta}}\, d{\theta} \bigg) \zeta^k.
$$ 
Since the function ${\theta}\mapsto \tilde
F(\xi,e^{i{\theta}})-\tilde h(e^{i{\theta}})$ is $C^\infty$ with
period $2{\pi}$, its Fourier series converges uniformly on $\R$ to the
function itself.  Hence, the sequence $(\tilde F_j)$ converges
uniformly on $\{\xi\}\times \T$ for each $\xi\in D_r$.  The
convergence is uniform on $D_t\times \T$ for each $t\in(1,r)$.
In fact, an integration by parts of the integral above shows that it  
can be estimated by
$$
\frac 1 {k^2}\max_{\xi\in D_t, {\theta}\in \R}\big| \partial^2
(\tilde F(\xi,e^{i{\theta}})-\tilde h(e^{i{\theta}}))
/\partial \theta^2 \big|,\qquad k\neq 0.
$$

Fixing $t\in (1,r)$, since $\tilde F(D_r\times \T)\subset Z$,
$F(\T\times\T)\subset X$, and $\tilde F_j \to\tilde F$ uniformly on
$D_t\times \T$, we have $\tilde F_j(D_t\times \T)\subset Z$ and
$\tilde F_j(\T\times\T)\subset \pi^{-1}(X)$ if $j$ is large enough.
For such $j$, define $F_j=\pi\circ\tilde F_j:D_t\times\T\to\P^n$.  The
$0$-th coordinate of $\tilde F$ is $\tilde
F_0(\xi,\cdot)=\|h-\gamma\|\xi/\varrho+1$, so the $0$-th coordinate of
$\tilde F_j$ is $\tilde F_{j0}(\xi,\cdot) =\chi_j\xi+1$, where
$\theta\mapsto \chi_j(e^{i\theta})$ is the $j$-th partial sum of the
Fourier series of $\theta\mapsto
\|h(e^{i\theta})-\gamma(e^{i\theta})\|/\varrho(e^{i\theta})$.  Hence,
$$
J(F_j(\cdot,\zeta))=\log|\chi_j(\zeta)|\to 
\log(\|h(\zeta)-\gamma(\zeta)\|/\varrho(\zeta))=J(F(\cdot,\zeta))
$$
uniformly for $\zeta\in \T$. Thus, by (3),
$$
\int_\T J(F_j(\cdot,\zeta))\, d\sigma(\zeta) 
<\int_\T v\circ h \, d\sigma+\epsilon
$$ 
for $j$ large enough.  We now fix $j$ so large that these
properties hold.

For every $\xi\in D_r$, the map $\zeta\mapsto \tilde
F_j(\xi,\zeta)-\tilde h(\zeta)$ has a pole of order at most $j$ at the
origin, and for every $\zeta\in D_r$, $\zeta\neq 0$, the map
$\xi\mapsto \tilde F_j(\xi,\zeta)-\tilde h(\zeta)$ has a zero at the
origin.  Hence, $(\xi,\zeta)\mapsto \tilde F_j(\xi\zeta^k,\zeta)$
extends to a holomorphic map $\overline \D\times \overline \D \to
\C^{n+1}$ for every $k\geq j$.

Since $\tilde F_j(0,\zeta)=\tilde h(\zeta)\in Z$ for all $\zeta\in
D_r$, $\zeta\neq 0$, there is $\delta >0$ such that $\tilde
F_j(\xi\zeta^k,\zeta)\in Z$ for all $k\geq j$ and $(\xi,\zeta)\in
D_\delta\times \overline \D$.  Since $\tilde F_j(\xi,\zeta)\in Z$ for
all $(\xi,\zeta)\in \overline \D\times \T$, there is $\tau<1$
such that $\tilde F_j(\xi,\zeta)\in Z$ for all $(\xi,\zeta)\in
\overline \D\times (\overline \D\setminus D_\tau)$, so $\tilde
F_j(\xi\zeta^k,\zeta)\in Z$ for all $(\xi,\zeta)\in \overline \D\times
( \overline \D\setminus D_\tau)$ and all $k\geq j$.  Choose $k\geq j$
large enough that $|\xi\zeta^k|<\delta$ for all $(\xi,\zeta)\in
\overline \D\times D_\tau$.  Then there is $s\in (1,t)$ such that
$\tilde F_j(\xi\zeta^k,\zeta)\in Z$ for all $(\xi,\zeta)\in D_s\times
D_s$.

Now define $\tilde G\in\O(D_s\times D_s,Z)$ by $\tilde G(\xi,\zeta)=
\tilde F_j(\xi\zeta^k,\zeta)$ and let $G=\pi\circ \tilde G$.  In the
proof of Theorem 1, we observed that if $\tilde f=(f_0,\dots,f_n)\in
\A_Z$ is a lifting of $f\in \A_{\P^n}$ and $f_0(0)\neq 0$, then
$$
J(f)=\int_\T \varphi\circ \tilde f\, d\sigma-\varphi(\tilde f(0)),
$$ 
where, as before, $\varphi(z)=\log|z_0|$ for $z\in\C^{n+1}$.
Now $\tilde G(0,\cdot)=\tilde h=(1,h)$, so $\varphi(\tilde
G(0,\cdot))=0$.  Therefore,
$$
\align
\int_\T J(G(\cdot,\zeta))\,d\sigma(\zeta)
&= \int_{\T^2} \varphi\circ \tilde G\, d(\sigma\times\sigma) 
=\dfrac 1{(2\pi)^2} \int_0^{2\pi}\int_0^{2\pi}
\varphi(\tilde F_j(e^{i(t+k\theta)},e^{i\theta}))\, dt\, d\theta \\
&= \int_{\T^2} \varphi\circ \tilde F_j\, d(\sigma\times\sigma) 
=\int_\T J(F_j(\cdot,\zeta))\, d\sigma(\zeta) 
<\int_\T v\circ h \, d\sigma+\epsilon.
\endalign
$$ 
By the Mean Value Theorem, there is ${\theta}_0\in [0,2{\pi}]$ such
that 
$$ \align
\int_\T J(G(\cdot,\zeta)\,d\sigma(\zeta)
&=\dfrac 1{(2\pi)^2}\int_0^{2{\pi}}\int_0^{2{\pi}}
{\varphi}(\tilde G(e^{i(\theta+t)},e^{it}))\, dt\, d{\theta} \\
&=\dfrac 1{2\pi}\int_0^{2{\pi}}
{\varphi}(\tilde G(e^{i(\theta_0+t)},e^{it})) \, dt.
\endalign $$ 
Now define $\tilde g(\zeta)=\tilde G(e^{i\theta_0}\zeta,\zeta)$ for
$\zeta\in D_s$, and $g=\pi\circ\tilde g$.  Then $\tilde g(0)=\tilde
G(0,0)=(1,h(0))$, so $g(0)=h(0)$, and $g(\T)\subset \pi(\tilde
G(\T\times\T))\subset X$, so $g\in\A_{\P^n}^X$.  Also,
$$\align
J(g)&=\int_\T\phi\circ\tilde g\,d\sigma
= \dfrac 1{2\pi}\int_0^{2{\pi}}
{\varphi}(\tilde G(e^{i{\theta}_0}e^{it},e^{it})) \, dt \\
&= \int_\T J(G(\cdot,\zeta))\, d\sigma(\zeta)
<\int_\T v\circ h \, d\sigma+\epsilon,
\endalign $$
and the proof is complete.
\qed\enddemo

\enddocument